\documentclass[12pt,final]{article}

\usepackage{amsmath,amssymb,amsthm}
\usepackage{color}
\usepackage[final]{graphicx}
\usepackage{bm} 
\usepackage{ascmac} 
\usepackage[capitalize]{cleveref} 
\usepackage{comment} 
\usepackage{stmaryrd} 
\usepackage{ulem} 
\usepackage{url}
\usepackage[color,draft]{showkeys} 

\newtheorem{thm}{Theorem}[section]
\newtheorem{crl}[thm]{Corollary}
\newtheorem{cnj}[thm]{Conjecture}
\newtheorem{prp}[thm]{Proposition}

\newtheorem{rmk}[thm]{Remark}
\newtheorem{dfn}[thm]{Definition}
\newtheorem{exm}[thm]{Example}

\newcommand{\st}{\stackrel}

\newcommand{\ra}{\rightarrow}
\newcommand{\mt}{\mapsto}

\makeatletter
\newcommand{\subjclass}[2][2020]{%
  \let\@oldtitle\@title%
  \gdef\@title{\@oldtitle\footnotetext{#1 \emph{Mathematics subject classification(s).} #2}}%
}
\newcommand{\keywords}[1]{%
  \let\@@oldtitle\@title%
  \gdef\@title{\@@oldtitle\footnotetext{\emph{Key words and phrases.} #1.}}%
}
\makeatother

\newcommand{\Tr}{\mathrm{Tr}\,}
\newcommand{\Nr}{\mathrm{N}_{n/n-1}\,}

\newcommand{\e}{\varepsilon_n}

\usepackage{here}

\title{Minimal relative units of the cyclotomic $\mathbb Z_2$-extension}
\author{
Tomokazu Kashio\thanks{Tokyo University of Science, \protect\url{kashio_tomokazu@ma.noda.tus.ac.jp}} \\
Hyuga Yoshizaki\thanks{Tokyo University of Science, \protect\url{yoshizaki.hyuga@gmail.com}}
}


\newtheorem*{cnj*}{Conjecture}
\newtheorem*{thm*}{Theorem}

\begin{document}

\maketitle

\begin{abstract}
Let $\mathbb B_n:=\mathbb Q(\cos(\pi/2^{n+1}))$. 
For the relative norm map $\Nr \colon \mathcal O_{\mathbb B_n}^\times \ra \mathcal O_{\mathbb B_{n-1}}^\times$ on the units group, 
we define $RE_n:=\mathrm{N}_{n/n-1}^{-1}(\{\pm 1\})$, $RE_n^+:=\mathrm{N}_{n/n-1}^{-1}(\{1\})$.
Komatsu conjectured that $\Tr \epsilon^2 \geq 2^n(2^{n+1}-1)$ for $\epsilon \in RE_n -\{\pm 1\}$.
Morisawa and Okazaki showed that it holds for $\epsilon \in RE_n -RE_n^+$.
In this paper we study the case $\epsilon \in RE_n^+$.
We conjecture that $\min \{\Tr \epsilon^2 \mid \epsilon \in RE_n^+-\{\pm 1\}\}= 2^n(1+8c_n)$, where $c_1:=2$ and $c_n:=2\cdot \mathrm{round}(2^n/5)$ ($n\geq 2$).
We show that this holds for $n\leq 6$ and that a ``half'' of this: $\min \{\Tr \epsilon^2 \mid \epsilon \in RE_n^+-\{\pm 1\}\} \leq 2^n(1+8c_n)$ holds for even $n$.
We also observe a relation to the class number problem.
\end{abstract}

\section{Introduction} \label{int}
Let $\mathbb B_n:=\mathbb Q(\cos(\pi/2^{n+1}))$, which is the $n$th layer of the cyclotomic $\mathbb Z_2$-extension over $\mathbb Q$.
We put $RE_n^+:=\operatorname{Ker} \Nr$, $RE_n^-:=\mathrm{N}_{n/n-1}^{-1}(\{-1\})$, $RE_n:=\mathrm{N}_{n/n-1}^{-1}(\{\pm 1\})=RE_n^+ \coprod RE_n^-$,  
where $\Nr \colon \mathcal O_{\mathbb B_n}^\times \ra \mathcal O_{\mathbb B_{n-1}}^\times$ denotes the relative norm map on the unit group.
Then Komatsu, in personal communication with Morisawa and Okazaki, stated the following conjecture.

\begin{cnj}[{\cite[Conjecture 1.1]{MO3}}] 
We have for $\epsilon \in RE_n -\{\pm 1\}$
\begin{align} \label{Kc}
\Tr \epsilon^2 \geq 2^n(2^{n+1}-1).
\end{align}
\end{cnj}

Morisawa and Okazaki showed that 

\begin{thm}[{\cite[Theorem 6.4]{MO3}}] \label{MO2}
{\rm Ineq.~(\ref{Kc})} holds for $\epsilon \in RE_n^-$.
\end{thm}

Namely the unsolved problem is Ineq.~(\ref{Kc}) for $\epsilon \in RE_n ^+-\{\pm 1\}$.
We provide the best possible refinement in this case as follows.
\begin{cnj*}[\cref{cnj}]
Let $c_1=2$, $c_n=2 \cdot \mathrm{round}\,(2^n/5)$ $(n\geq 2)$ where $\mathrm{round}\,(x)$ denotes the nearest integer to $x$.
Then we have
\begin{align} \label{refinement}
\min \{\Tr \epsilon^2 \mid \pm 1 \neq \epsilon \in RE_n^+\}=2^n(1+8c_n).
\end{align}
\end{cnj*}
The first few terms of $c_n$ are $c_1=2$, $c_2=2$, $c_3=4$, $c_4=6$, $c_5=12,\dots$.
We also present some partial results.
\begin{thm*}[\cref{lb}]
For $n=1,3,5$ or for even $n$, there exists $u_n \in RE_n^+-\{\pm 1\}$ satisfying 
\begin{align*}
\Tr u_n^2 =2^n(1+8c_n).
\end{align*}
\end{thm*}
Hence a ``half'' of \rm\cref{refinement}: $\min \{\Tr \epsilon^2 \mid \pm 1 \neq \epsilon \in RE_n^+\}\leq 2^n(1+8c_n)$ holds for such $n$.
\begin{thm*}[\cref{num}]
{\rm\cref{refinement}} holds for $n\leq 6$. 
\end{thm*}
The proof of \cref{num} relies on the fact that the class number $h_n$ of $\mathbb B_n$ is $1$.
On the other hand, in \cref{alg}, we also provide a proof for $n\leq 3$ without using any information of $h_n$.

We also see a relation between our Conjecture and the class numbers in \S \ref{h3=1}, \ref{lind}.
Weber's class number problem asks whether $h_n=1$ for all $n$ and some partial results follows by studying the unit group.
For example, by using \cref{MO2} concerning $RE_n^-$, Fukuda and Komatsu \cite[Theorem 1.3]{FK3} showed that 
\begin{align} \label{fkf}
\text{$l \nmid h_n$ for all $n$ and for all primes $l$ with $l\not \equiv \pm 1 \bmod 32$.}
\end{align}
We may observe a ``similar'' phenomena also for $RE_n^+$.
Morisawa and Okazaki \cite[Proposition 6.6]{MO3} showed that
\begin{align} \label{mof}
\min \{\Tr \epsilon^2 \mid \pm 1 \neq \epsilon \in RE_n^+\}\geq 2^n\cdot 17 \quad (n\geq 2).
\end{align}
The second author \cite[Remark in \S 5.1]{Yo} showed that
\begin{center}
Ineq.~(\ref{mof}) implies $h_2/h_1=1$.
\end{center}
We generalize these results as follows. 

\begin{thm*}[\cref{Y}]
We have 
\begin{align} \label{yf}
\min \{\Tr \epsilon^2 \mid \pm 1 \neq \epsilon \in RE_n^+\}\geq 2^n \cdot 33  \quad (n\geq 3).
\end{align}
\end{thm*}

\begin{thm*}[\cref{n=3}]
{\rm Ineq.~(\ref{yf})} implies $h_3/h_2=1$.
\end{thm*}

In \S \ref{lind} we introduce some numerical results:
\begin{enumerate}
\item When $n=4,5$, \cref{refinement} implies the $l$-indivisibility of $h_n/h_{n-1}$ for several primes $l$ (\S \ref{ind}).
\item When $n=7$, \cref{refinement} implies the $l$-indivisibility of $h_7/h_{7-1}$ for the first $1000$ primes $l$ satisfying $l>10^9$, $l\equiv 65 \bmod 128$ (\S \ref{new}).
\end{enumerate}
The known results for the $l$-indivisibility is as follows.
\begin{align} \label{previous}
\text{if $n<7$ or $l\not\equiv \pm 1 \bmod 64$ or $l\leq10^9$, then a prime $l$ does not divide $h_n$}.
\end{align}
We note that the primes $l$ in the case (ii) are out of this range.

\section{Minimal relative units} \label{min}

Let $\mathbb B_n$ be the $n$th layer of the cyclotomic $\mathbb Z_2$-extension over $\mathbb Q$. More explicitly we have
\begin{align*}
\mathbb B_n=\mathbb Q(X_n), \quad X_n:=2\cos\left(\frac{2\pi}{2^{n+2}}\right).
\end{align*}
In this paper, we fix a generator $\sigma$ of $G_n:=\mathrm{Gal}(\mathbb B_n/\mathbb Q)\cong \mathbb Z/2^n\mathbb Z$ by
\begin{align*}
\sigma \colon 2\cos\left(\frac{2\pi}{2^{n+2}}\right) \mt 2\cos\left(\frac{3 \cdot 2\pi}{2^{n+2}}\right).
\end{align*}

\begin{dfn}
Let $E_n$ be the unit group of $\mathbb B_n$.
We consider the following subgroups:
\begin{align*}
RE_n^+&:=\{\epsilon \in E_n \mid \Nr \epsilon =1 \}, \\
RE_n&:=\{\epsilon \in E_n \mid \Nr \epsilon =\pm 1 \}, \\
A_n&:=\langle \pm 1, \e \rangle_{\mathbb Z[G_n]}=\left\{\pm \prod_{i=0}^{2^{n-1}-1} \sigma^i(\e)^{m_i}\mid  m_i \in \mathbb Z\right\} \quad \text{for}\quad \e:=\frac{X_n+1}{X_n-1}.
\end{align*}
Here $\Nr\colon \mathbb B_n \ra \mathbb B_{n-1}$ denotes the relative norm map. 
\end{dfn}
We have $A_n \subset RE_n^+$ since $\Nr \sigma^i(\e) =\sigma^i(\frac{X_n+1}{X_n-1} \cdot \frac{-X_n+1}{-X_n-1})=1$.
We embed $\mathbb B_n$ into $\mathbb R^{2^n}$ as usual:
\begin{align*}
\mathbb B_n \ra \mathbb R^{2^n}, \quad x \mt (\sigma^i(x))_{0\leq i \leq 2^n-1}.
\end{align*}
In particular, $\sqrt{\Tr x^2}$ is equal to the length of $x$ in $\mathbb R^{2^n}$.
The ring of integers $\mathcal O_{\mathbb B_n} =\mathbb Z[X_n]$ has an orthogonal basis $\{b_i \mid 0\leq i \leq 2^n-1\}$:
\begin{align} \label{ob}
b_i:= 
\begin{cases}
1 & (i=0) \\
2\cos\left(\frac{i*2\pi}{2^{n+2}}\right) & (1\leq i \leq 2^n-1)
\end{cases}, \quad 
\mathrm{Tr}\, (b_ib_j)=
\begin{cases}
0 & (i\neq j) \\
2^n & (i=j=0) \\
2^{n+1} & (i=j>0)
\end{cases}.
\end{align}
In this paper, we repeatedly use the following relations:
\begin{align*}
b_0 b_i=b_i, \quad b_i b_j=b_{i+j}+b_{i-j}, \quad b_i^2=2+b_{2i} \qquad (1\leq i,j \leq 2^n-1, \ i\neq j),
\end{align*}
where we regard that 
\begin{align*}
b_{2^n}=0, \quad b_{-k}:=b_k, \quad b_{2^n+k}:=-b_{2^n-k} \quad (1 \leq k \leq 2^n-1).
\end{align*}

The following conjecture and the partial results below are the main results in this paper.

\begin{cnj} \label{cnj}
We define $c_n$ for $n \in \mathbb N$ by 
\begin{align*}
c_1&:=2, \\
c_n&:=2 \cdot \mathrm{round}\,(2^n/5)=
\begin{cases}
2 (2^n-1)/5 & (n\equiv  0\mod 4) \\
2 (2^n-2)/5 & (n\equiv  1\mod 4) \\
2 (2^n+1)/5 & (n\equiv  2\mod 4) \\
2 (2^n+2)/5 & (n\equiv  3\mod 4) \\
\end{cases} \qquad (n\geq 2).
\end{align*}
Here $\mathrm{round}\,(x)$ denotes the nearest integer to $x$.
Then we have
\begin{align*}
\min \{\Tr \epsilon^2 \mid \pm 1 \neq \epsilon \in RE_n^+\}=2^n(1+8c_n).
\end{align*}
\end{cnj}
For example, $c_1=2$, $c_2= 2$, $c_3=4$, $c_4=6$, $c_5=12$, $c_6=26$, $c_7= 52$, $c_8= 102$, $c_9= 204$, $c_{10}= 410$.
Hereinafter in this section, we present partial results (Theorems \ref{Y}, \ref{lb} and \cref{alg}) for \cref{cnj}.
First, we generalize Morisawa-Okazaki's Ineq.~(\ref{mof}) a little. 

\begin{thm} \label{Y} 
We have for $n\geq 3$
\begin{align*}
\min \{\Tr \epsilon^2 \mid \pm 1 \neq \epsilon \in RE_n^+\} \geq 2^n \cdot 33.
\end{align*}
\end{thm}

\begin{proof}
Let $\epsilon \in RE_n^+$, $\neq \pm 1$. Write
\begin{align*}
\epsilon=\sum_{i=0}^{2^n-1}a_ib_i \quad (a_i \in \mathbb Z).
\end{align*}
We have by $\epsilon \in RE_n^+$
\begin{align} \label{e-o}
\Nr \epsilon =\biggl(\sum_{2\mid i}a_ib_i\biggr)^2-\biggl(\sum_{2\nmid i}a_ib_i\biggr)^2=1.
\end{align}
\cite[Lemma 6.2]{MO3} states that 
\begin{align} \label{oe}
\text{$a_0 $ is odd, $a_i$ ($i\neq 0$) are even}.
\end{align}
We claim that it suffices to show that 
\begin{itemize}
\item[(a)] at least four $a_i$'s are not equal to $0$ for non-zero $i$, or, 
\item[(b)] at least two $a_i$'s are not equal to $0$ for odd $i$.
\end{itemize}
First we note that 
\begin{align*} 
\Tr \epsilon^2=2^n(a_0^2+2a_1^2+\cdots +2a_{2^n-1}^2).
\end{align*}
by (\ref{ob}).
The statement (a) implies the assertion since we have 
\begin{align*}
2^n(a_0^2+2a_1^2+\cdots +2a_{2^n-1}^2)\geq 2^n(1+2 \cdot 4 \cdot 2^2)=2^n\cdot 33
\end{align*}
by (\ref{oe}). Now assume (b). By taking the trace of (\ref{e-o}), we have
\begin{align*}
2^na_0^2+2^{n+1}\sum_{2\mid i \neq 0}a_i^2-2^{n+1} \sum_{2\nmid i}a_i^2=2^n.
\end{align*}
It follows that  
\begin{align*}
2^n(a_0^2+2a_1^2+\cdots +2a_{2^n-1}^2)=2^n+2^{n+2}\sum_{2\nmid i}a_i^2.
\end{align*}
This is greater than or equal to $2^n+2^{n+2}\cdot 2\cdot 4=2^n\cdot 33$ by (\ref{oe}) and (b) as desired.

Recall that $\pm 1 \neq \epsilon \in RE_n^+$. In particular $\epsilon  \in \mathbb B_n-\mathbb B_{n-1}$, so at least one $a_i$ is not equal to $0$ for odd $i$.
We may assume $i=1$ by considering the Galois action.
If there exists at least one more odd $i$ satisfying $a_i\neq 0$, then (b) holds.
Assume that $a_i=0$ for odd $i\neq 1$. Then (\ref{e-o}) becomes 
\begin{align*}
\biggl(\sum_{2\mid i}a_ib_i\biggr)^2=(1+2a_1^2)+a_1^2b_2.
\end{align*}
By (\ref{oe}), we have $a_0 \neq 0$. There exists at least one more even $i_1$ satisfying $a_{i_1}\neq 0$, 
because otherwise it follows that $a_1^2b_2=a_0^2-1-2a_1^2 \in \mathbb Z$.
This is a contradiction for $b_2 \in \mathbb B_{n-1}-\mathbb B_{n-2}$ and $n\geq 3$.
Once again, we see that there exists at least one more even $i_2$ satisfying $a_i\neq 0$, because otherwise it follows that 
\begin{align*}
(a_0^2+2a_{i_1}^2)+2a_0a_{i_1}b_{i_1}+a_{i_1}^2b_{2i_1}=(1+2a_1^2)+a_1^2b_2.
\end{align*}
Then we have ``$b_2=b_{i_1}=- b_{2i_1}$'' or ``$b_2=b_{i_1}$, $b_{2i_1} \in \mathbb Z$'', that is, ``$i_1=2$, $2i_1=2^{n+1}-2$'' or ``$i_1=2$, $2i_1=2^n$'', which is a contradiction for $n\geq 3$.

Now $\epsilon$ has at least three non-zero coefficients $a_1,a_{i_1},a_{i_2}$ with $2\mid i_1,i_2$, other than $a_0$.
We assume for the contradiction that these are all non-zero ones.
In particular (\ref{e-o}) becomes  
\begin{align}
&(a_0^2+2a_{i_1}^2+2a_{i_2}^2)+2a_0a_{i_1}b_{i_1}+2a_0a_{i_2}b_{i_2}+a_{i_1}^2b_{2i_1}+a_{i_2}^2b_{2i_2}+2a_{i_1}a_{i_2}b_{i_1+i_2}+2a_{i_1}a_{i_2}b_{i_1-i_2} \notag \\
&=(1+2a_1^2)+a_1^2b_2. \label{e-o2}
\end{align}
We consider three cases: ``$i_1\equiv i_2 \equiv 0 \bmod 4$'', ``$i_1\equiv i_2 \equiv 2 \bmod 4$'', ``$i_1\equiv 2 \bmod 4$, $i_2 \equiv 0 \bmod 4$''.
First assume that $i_1\equiv i_2 \equiv 0 \bmod 4$.
Then we have $2i_1,2i_2, i_1\pm i_2 \equiv 0 \bmod 4$. Therefore there does not exist any term in the left-hand side corresponding to $a_1^2b_2$ in the right-hand side, which is a contradiction. 
Next assume that $i_1\equiv i_2 \equiv 2 \bmod 4$.
We have $2i_1,2i_2, i_1\pm i_2 \equiv 0 \bmod 4$. 
Therefore the relation (\ref{e-o2}) implies 
\begin{align*}
2a_0a_{i_1}b_{i_1}+2a_0a_{i_2}b_{i_2}=a_1^2b_2.
\end{align*}
This follows, for example, by considering the quotient vector space $\mathbb B_{n-1}/\mathbb B_{n-2}$.
Then we have $i_1=i_2=2$, which is a contradiction.
Finally assume that $i_1\equiv 2 \bmod 4$, $i_2 \equiv 0 \bmod 4$.
We have $i_1\pm i_2 \equiv 2 \bmod 4$, $2i_1,2i_2 \equiv 0 \bmod 4$. Similarly as above we obtain  
\begin{align*}
2a_0a_{i_1}b_{i_1}+2a_{i_1}a_{i_2}b_{i_1+i_2}+2a_{i_1}a_{i_2}b_{i_1-i_2}=a_1^2b_2. 
\end{align*}
We have $i_1+i_2 \notin \{\pm 2,\pm(2^{n+1}-2)\}$ by $2\leq i_1 \leq 2^n-2$, $4\leq i_2 \leq 2^n-4$. That is, $|b_2| \neq |b_{i_1+i_2}|$.
Then there are two possible cases: 
\begin{align*}
|b_2|=|b_{i_1}| \neq |b_{i_1+i_2}|=|b_{i_1-i_2}| \quad \text{or}  \quad |b_2|=|b_{i_1-i_2}| \neq |b_{i_1+i_2}|=|b_{i_1}|.
\end{align*}
If the former one holds, then we have $i_1=2$ and $i_1+i_2 =-(i_1-i_2)$, which is a contradiction.
If the latter one holds, we have 
\begin{align*}
2a_0a_{i_1}b_{i_1}+2a_{i_1}a_{i_2}b_{i_1+i_2}=0,
\end{align*}
which implies $|a_0|=|a_{i_2}|$. This is a contradiction for (\ref{oe}). Then the assertion is clear. 
\end{proof}

\begin{rmk}
\begin{enumerate}
\item The above proof is independent of any information of the class number $h_n$ of $\mathbb B_n$.
Oppositely, we show that {\rm\cref{Y}} implies $h_3/h_2=1$ in {\rm\S \ref{h3=1}}.
\item The strategy of the above proof is counting the number of non-zero coefficients $a_i$ of a relative unit $\epsilon = \sum_{i=0}^{2^n-1} a_ib_i \in RE_n^+$ by a combinatorial argument,
and showing that the number is greater than or equal to $c_k$ if $n\geq k$, for $k=3$.
The same proof works for $k=4$, although we used a computer.
\end{enumerate}
\end{rmk}

For small $n$ or even $n$, we obtain (a candidate of) the minimal unit $\in RE_n^+$ explicitly.

\begin{thm} \label{lb} For $n=1,3,5$, we put 
\begin{align*}
u_1&:=\varepsilon_1=3b_0+2b_1, \\
u_3&:=\varepsilon_3 \sigma(\varepsilon_3)=b_0+ 2(b_1+b_2+b_5+b_6), \\
u_5&:=\varepsilon_5 \sigma^2(\varepsilon_5)\\
&=b_0+2(b_{11}+b_{12}-b_{14}-b_{15}+b_{17}+ b_{18}+ b_{19}+ b_{20}-b_{22}-b_{23}+b_{25}+ b_{26}).
\end{align*}
For $n \in 2\mathbb N$, we put
\begin{align*}
u_n:=b_0+ (-1)^{\frac{n}{2}}2\sum_{i=\lceil \frac{2^{n+1}}{5}\rceil}^{\lfloor \frac{2^{n+2}}{5}\rfloor} b_i.
\end{align*}
Here $\lceil \ \rceil$, $\lfloor \ \rfloor$ denote the ceiling function, the floor function, respectively. 
Then we have 
\begin{align*}
\Tr u_n^2=2^n(1+8c_n) \quad (n=1,3,5 \text{ or } n \in 2\mathbb N). 
\end{align*}
\end{thm}
Hence a ``half'' of \cref{cnj} holds for such $n$:
\begin{align*}
\min \{\Tr \epsilon^2 \mid \pm 1 \neq \epsilon \in RE_n^+\}\leq 2^n(1+8c_n) \quad (n=1,3,5 \text{ or } n \in 2\mathbb N).
\end{align*}

\begin{proof}
The cases $n=1,3,5$ follow from a direct calculation, by noting that (\ref{ob}) implies
\begin{align} \label{tr}
\Tr \biggl(\sum_{i=0}^{2^n-1} c_i b_i \biggr)^2=2^n\biggl(c_0^2+2\sum_{i=1}^{2^n-1} c_i^2\biggr)  \quad (c_i \in \mathbb Z).
\end{align}

For even $n$, easily see that 
\begin{align} \label{f-c}
c_n=\lfloor \tfrac{2^{n+2}}{5}\rfloor-\lceil \tfrac{2^{n+1}}{5}\rceil+1.
\end{align}
It follows that $\Tr u_n^2=2^n(1+8c_n)$ by (\ref{tr}).
Hence it suffices to show that $\Nr u_n=1$.
Let $s:=\lceil \frac{2^{n+1}}{5}\rceil$, $t:=\lfloor \frac{2^{n+2}}{5}\rfloor$, $b(n):=b_n$.
We can write
\begin{align*}
\mathrm{N}_{n/n-1}\,u_n-1 &=\left(1+ (-1)^{\frac{n}{2}}2\sum_{i=s}^{t} b(i)\right)\left(1+ (-1)^{\frac{n}{2}}2\sum_{i=s}^{t} (-1)^i b(i)\right) -1\\
&=(-1)^{\frac{n}{2}}4\sum_{s \leq 2k \leq t} b(2k)  + 4\sum_{s \leq 2k \leq t}b(2k)^2 + 8\sum_{s \leq 2k<2l \leq t} b(2k)b(2l) \\
&\quad -4\sum_{s \leq 2k+1 \leq t}b(2k+1)^2 -8\sum_{s \leq 2k+1<2l+1 \leq t} b(2k+1)b(2l+1).
\end{align*}
The sum of the second and forth terms in the most right-hand side is equal to
\begin{align*}
4\sum_{s \leq 2k \leq t}(b(4k)+2)  -4\sum_{s \leq 2k+1 \leq t}(b(4k+2)+2)
=4\sum_{s \leq 2k \leq t}b(4k)  -4\sum_{s \leq 2k+1 \leq t}b(4k+2).
\end{align*}
since (\ref{f-c}) implies that the parities of $s,t$ are even-odd or odd-even.
The sum of the third and fifth terms is equal to 
\begin{align*}
&8\sum_{s \leq 2k<2l \leq t} (b(2k+2l)+b(2k-2l)) -8\sum_{s \leq 2k+1<2l+1 \leq t} (b(2k+2l+2)+b(2k-2l)) \\
&=8\sum_{s \leq 2k<2l \leq t} b(2k+2l) -8\sum_{s \leq 2k+1<2l+1 \leq t} b(2k+2l+2),
\end{align*}
by the parities of $s,t$ again.
Hence it suffices to show that 
\begin{align} 
&(-1)^{\frac{n}{2}} \sum_{s \leq 2k \leq t} b(2k)  +\sum_{s \leq 2k \leq t}b(4k) -\sum_{s \leq 2k+1 \leq t}b(4k+2) \notag \\ 
&\quad + 2\sum_{s \leq 2k<2l \leq t} b(2k+2l) -2\sum_{s \leq 2k+1<2l+1 \leq t} b(2k+2l+2) \label{0} 
\end{align}
is equal to $0$. We divide it into two cases.
First assume that $n \equiv 2 \mod 4$. 
Then $s$ is even and $t$ is odd. Therefore $s \leq 2k+1  \leq t$ is equivalent to $s \leq 2k \leq t-1$, 
and $s \leq 2k+1<2l+1  \leq t$ is equivalent to $s \leq 2k<2l \leq t-1$, respectively.
Then (\ref{0}) becomes  
\begin{align*}
&- \sum_{s \leq 2k \leq t-1} b(2k)  +\sum_{s \leq 2k \leq t-1}b(4k) -\sum_{s \leq 2k \leq t-1}b(4k+2) \\
&\quad + 2\sum_{s \leq 2k<2l \leq t-1} b(2k+2l) -2\sum_{s \leq 2k<2l \leq t-1} b(2k+2l+2).
\end{align*}
Since we have  
\begin{align*}
\sum_{s \leq 2k<2l \leq t-1} b(2k+2l) -\sum_{s \leq 2k<2l \leq t-1} b(2k+2l+2) &=\sum_{s\leq 2k \leq t-1} b(4k+2)  - \sum_{s+t+1\leq 2k \leq 2t}  b(2k), \\
\sum_{s \leq 2k \leq t-1}b(4k) +\sum_{s \leq 2k \leq t-1}b(4k+2) &= \sum_{2s \leq 2k \leq 2t} b(2k),
\end{align*}
the problem is reduced to showing that  
\begin{align*}
-\sum_{s \leq 2k \leq t-1} b(2k)  + \sum_{2s \leq 2k \leq 2t}b(2k) - 2\sum_{s+t+1\leq 2k \leq 2t}  b(2k) =0.
\end{align*}
Let $c:=\frac{2^{n+2}}{10}$. For even $n$ (not only for $n\equiv 2 \bmod 4$), we see that
\begin{align}
&\bullet \text{ $s$ is the least integer $\geq  c $, $t$ is the greatest integer $\leq 2c $}, \notag \\[3pt]
&\bullet \text{ $2s$ is the least even integer $\geq 2 c $, $2t$ is the greatest even integer $\leq 4 c $}, \label{gili} \\[3pt]
&\bullet \text{ $s+t+1$ is the least even integer $\geq 3 c $}. \notag
\end{align}
Therefore the left-hand side becomes
\begin{align*}
&-\sum_{c  \leq 2k \leq 2c } b(2k)  + \sum_{2 c  \leq 2k \leq 4 c }b(2k) 
- 2\sum_{3 c \leq 2k \leq 4 c }  b(2k) \\
&=-\sum_{c  \leq 2k \leq 2c } b(2k)  + \sum_{2 c  \leq 2k \leq 3 c }b(2k) - \sum_{3 c \leq 2k \leq 4 c }  b(2k) 
\end{align*}
The last sum is equal to $0$ since we have
\begin{align*}
&\sum_{c  \leq 2k \leq 2c } b(2k) +\sum_{3 c \leq 2k \leq 4 c }  b(2k) =0, \\
&\sum_{2 c  \leq 2k \leq 3 c }b(2k)=\sum_{2 c  \leq 2k < 2^n }b(2k)+b(2^n)+\sum_{2^n < 2k < 3c }b(2k)=0
\end{align*}
by $b(2^n+k)=-b(2^n-k)$ and $b(2^n)=0$.

Next assume that $n \equiv 0 \mod 4$, which implies $s$ is odd and $t$ is even.
Then (\ref{0}) becomes 
\begin{align}
&\sum_{s+1 \leq 2k \leq t} b(2k)  +\sum_{s+1 \leq 2k \leq t}b(4k) -\sum_{s-1 \leq 2k \leq t-2}b(4k+2) \notag \\
&\quad + 2\sum_{s+1 \leq 2k<2l \leq t} b(2k+2l) -2\sum_{s-1 \leq 2k<2l \leq t-2} b(2k+2l+2). \label{00}
\end{align}
In this case we have
\begin{align*}
\sum_{s+1 \leq 2k<2l \leq t} b(2k+2l) -\sum_{s-1 \leq 2k<2l \leq t-2} &b(2k+2l+2) \\
&=\sum_{s-1 \leq 2k \leq t-2} b(4k+2) -\sum_{2s \leq 2k \leq s+t-1} b(2k), \\
\sum_{s+1 \leq 2k \leq t}b(4k) +\sum_{s-1 \leq 2k \leq t-2}b(4k+2)&=\sum_{2s \leq 2k \leq 2t}b(2k).
\end{align*}
Hence (\ref{00}) is equal to 
\begin{align} \label{000}
\sum_{s+1 \leq 2k \leq t} b(2k)  +\sum_{2s \leq 2k \leq 2t}b(2k)  -2\sum_{2s \leq 2k \leq s+t-1} b(2k).
\end{align}
By (\ref{gili}), we can rewrite (\ref{000}) as 
\begin{align*}
&\sum_{ c  \leq 2k \leq 2 c } b(2k)  +\sum_{2 c  \leq 2k \leq 4 c }b(2k)  
-2\sum_{2 c  \leq 2k \leq 3 c } b(2k) \\
&=\sum_{ c  \leq 2k \leq 2 c } b(2k)  +\sum_{3 c  \leq 2k \leq 4 c }b(2k)  
-\sum_{2 c  \leq 2k \leq 3 c } b(2k),
\end{align*}
which is equal to $0$ by $b(2^n+k)=-b(2^n-k)$ and $b(2^n)=0$.
Then the assertion is clear.
\end{proof}

We obtain the following corollary by Ineq.~(\ref{mof}), Theorems \ref{Y}, \ref{lb} (and a trivial argument for $n=1$).
\begin{crl} \label{alg}
{\rm\cref{cnj}} holds true for $n=1,2,3$.
\end{crl}

\cref{cnj} should be proved without studying the class number $h_n$ of $\mathbb B_n$, as we seen above.
On the other hand, we have 
\begin{align} \label{kn}
k_n:=\frac{h_n}{h_{n-1}}=[RE_n^+:A_n].
\end{align}
This follows from, for example, \cite[Theorem 8.2, Proposition 8.11]{Wa}, \cite[(1), (4)]{H2}. For a proof, see \cite[\S 4.1]{Yo}. 
Besides, we have $h_n=1$ for $n\leq 6$, so $RE_n^+=A_n$ for the same $n$. 
Since $A_n$ is given explicitly, we can verify \cref{cnj} numerically for such $n$, as follows.

Assume that $u \in RE_n^+$ satisfies 
\begin{align} \label{tru}
\Tr u^2 \leq 2^n(1+8c_n).
\end{align}
We put
\begin{align*}
x_i:=\log |\sigma^i(u)| \in \mathbb R \quad (0\leq i \leq 2^{n-1}-1).
\end{align*}
Since $\Nr \tau(u)= \tau(u)\sigma^{2^{n-1}}(\tau(u))=1$ for $\tau \in G_n$, the inequality (\ref{tru}) turns into  
\begin{align*}
\sum_{i=0}^{2^{n-1}-1} (e^{2x_i} +e^{-2x_i}) \leq 2^n(1+8c_n).
\end{align*}
We consider the logarithmic embedding 
\begin{align*}
RE_n^+/\{\pm 1\} \hookrightarrow \mathbb R^{2^{n-1}}, \quad \epsilon \mt (\log(|\sigma^i(\epsilon)|))_{i=0,1,\dots,2^{n-1}-1}.
\end{align*}
Then the square of the length of the image of $u$ is given by 
\begin{align*}
\sum_{i=0}^{2^{n-1}-1} (\log |\sigma^i(u)|)^2=\sum_{i=0}^{2^{n-1}-1} x_i^2.
\end{align*}
We put
\begin{align*}
L_n:=\max\left \{\sum_{i=0}^{2^{n-1}-1} x_i^2 \, \middle | \, x_i \in \mathbb R,\ \sum_{i=0}^{2^{n-1}-1} (e^{2x_i} +e^{-2x_i}) \leq 2^n(1+8c_n) \right\}.
\end{align*}
Namely, the condition (\ref{tru}) implies 
\begin{align} \label{trlu}
\sum_{i=0}^{2^{n-1}-1} (\log |\sigma^i(u)|)^2 \leq L_n.
\end{align}
Now we assume that $RE_n^+=A_n$, which is equivalent to $k_n:=\frac{h_n}{h_{n-1}}=1$ by (\ref{kn}) ($n\leq 6$ is a sufficient condition). 
Then we may write
\begin{align*}
u=\prod_{j=0}^{2^{n-1}-1} \sigma^j(\e)^{n_j} \quad (n_j \in \mathbb Z).
\end{align*}
Therefore (\ref{trlu}) is equivalent to 
\begin{align} 
&M[\mathbf n]:={}^t\mathbf nM \mathbf n \leq L_n,  \label{cnd} \\
&M:=\left[ \sum_{k=0}^{2^{n-1}-1} \log |\sigma^{k+i}(\e)| \log |\sigma^{k+j}(\e)| \right]_{i,j=0,1,\dots,2^{n-1}-1}, \notag \\
&\mathbf n:=[n_i]_{i=0,1,\dots,2^{n-1}-1}. \notag
\end{align}
We can find all such vectors $\mathbf n$ by the Fincke-Pohst algorithm (actually, we used the command {\tt qfminim} of PARI/GP).
Here the value of $L_n$ is given as follows:  
Assume that $x_i$ satisfies $\sum_{i=0}^{2^{n-1}-1} (e^{2x_i} +e^{-2x_i})=2^na$ for a fixed $a$. Note that $a \geq 1$ since $y+y^{-1}\geq 2$ for $y \in \mathbb R$.
Then the Lagrange multiplier theorem says that the function $\sum_{i=0}^{2^{n-1}-1} x_i^2$ takes the maximum value 
only when $(x_i)_i=\lambda (e^{2x_i} -e^{-2x_i})_i$ for some $\lambda \in \mathbb R$.
The solutions of $x=\lambda(e^{2x} -e^{-2x})$ are of the form of $x=\pm b$ with $b \geq 0$ 
since $e^{2x} -e^{-2x}$ ($x \geq 0$) is a convex function  and $x,e^{2x} -e^{-2x}$ are odd functions.
It follows that $e^{2x_i} +e^{-2x_i}$ is constant for all $i$, that is, $e^{2x_i} +e^{-2x_i}=2a$.
Namely, $\sum_{i=0}^{2^{n-1}-1} x_i^2$ takes the maximum value when 
\begin{align*}
x_i=\pm \frac{\log(a-\sqrt{a^2-1})}{2}.
\end{align*}
Therefore we see that 
\begin{align*}
L_n&=\max \left\{2^{n-3} \left(\log\left(a-\sqrt{a^2-1}\right)\right)^2 \, \middle | \,  1 \leq a \leq 1+8c_n \right\}\\
&=2^{n-3} \left(\log\left(1+8c_n-\sqrt{16c_n+64c_n^2}\right)\right)^2 
\end{align*}
In fact, we have $L_1=3.107\ldots$, $L_2=6.214\ldots$, $L_3=17.55\ldots$, $L_4=42.04\ldots$, $L_5=111.0\ldots$, $L_6=291.4\ldots$, $L_7=723.8\ldots$. 

When $n\leq 6$, we confirmed that $u$ does not satisfies $\Tr u^2 < 2^n(1+8c_n)$ for any $\mathbf n \neq \mathbf 0$ satisfying (\ref{cnd}):
for example, let $n=6$. 
Then the number of vectors $\mathbf n\neq \mathbf 0$ satisfying (\ref{cnd}) is $290624$. 
We computed $\Tr (\prod_{i=0}^{2^{5}-1} \sigma(\e)^{n_i})^2$ for such $\mathbf n$ and checked that the minimal value is equal to $2^6(1+8c_6)$.
To summarize, by numerical computation and by using $k_n=1$, we have the following.

\begin{thm} \label{num}
{\rm\cref{cnj}} holds true for $n\leq 6$.
\end{thm}

\begin{rmk}
\begin{enumerate}
\item When $n\leq 6$, all $\epsilon \in RE_n^+$ satisfying $\Tr \epsilon^2 =2^n(1+8c_n)$ are the conjugates of $u_n$ given in {\rm\cref{lb}}.
\item We can not confirm the case $n>6$ due to the limit of computer power. 
\end{enumerate}
\end{rmk}

\section{Relation to $h_n=1$ when $n\leq 3$} \label{h3=1}

There are many partial results supporting Weber's class number problem obtained by studying the unit group.
More directly, the second author proved the following.
We put $k_n:=\frac{h_n}{h_{n-1}}$, where $h_n$ denotes the class number of $\mathbb B_n$.

\begin{thm}[{\cite[Remark in \S 3.3]{Yo}}]
{\rm Ineq.~(\ref{mof})} implies $k_2=1$.
\end{thm}

In this section, we generalize this result to $n=3$ as follows.

\begin{thm} \label{n=3}
{\rm\cref{Y}} implies $k_3=1$.
\end{thm}

\begin{proof}
Assume for contradiction that there exists $u \in RE^+_3-A_3$. 
Since $[RE^+_3:A_3]<\infty$, we can write 
\begin{align*}
u:=\prod_{j=0}^{3} \sigma^j(\varepsilon_3)^{x_j}, \quad (x_0,x_1,x_2,x_3) \in \mathbb Q^4-\mathbb Z^4.
\end{align*}
We may replace $(x_0,x_1,x_2,x_3)$ with $(x_0',x_1',x_2',x_3')$ so that $x_i \equiv x_i' \bmod \mathbb Z$.
Therefore, by putting  
\begin{align*}
&T(x_0,x_1,x_2,x_3):=\sum_{i=0}^7 \prod_{j=0}^{3} |\sigma^{i+j}(\varepsilon_3)|^{2x_j}, \\ 
&L:=\min_{\alpha_0 \in \mathbb R} \max_{\alpha_0 \leq x_0 \leq \alpha_0+1 } 
\min_{\alpha_1 \in \mathbb R} \max_{\alpha_1 \leq x_1 \leq \alpha_1+1 } 
\min_{\alpha_2 \in \mathbb R} \max_{\alpha_2 \leq x_2 \leq \alpha_2+1 } 
\min_{\alpha_3 \in \mathbb R} \max_{\alpha_3 \leq x_3 \leq \alpha_3+1 } T(x_0,x_1,x_2,x_3),
\end{align*}
it suffices to show that 
\begin{align} \label{bd3}
L  < 2^3(1+8c_3)=264
\end{align}
since we have 
\begin{align*}
T(x_0,x_1,x_2,x_3)=\Tr u^2\geq 264
\end{align*}
for $\pm 1 \neq u \in RE_3^+$ by \cref{Y}.

First we show that $T(x_0,x_1,x_2,x_3)$ is a convex function.
In particular, a set 
\begin{align*}
S:=\{(x_0,x_1,x_2,x_3,y) \mid y\geq T(x_0,x_1,x_2,x_3)\}
\end{align*}
is convex.
We can reduce it to the convexity of a function of the form $a^{x_0}b^{x_1}c^{x_2}d^{x_3}$ since the sum of convex functions is again convex.
Its Hessian matrix is equal to 
\begin{align*}
4a^{2 x_0} b^{2 x_1} c^{2 x_2} d^{2 x_3}
\begin{pmatrix}
(\log a) ^2&
\log a  \log b &
\log a  \log c &
\log a  \log d\\ 
\log a  \log b &
(\log b) ^2&
\log b  \log c &
\log b  \log d\\ 
\log a  \log c &
\log b  \log c &
(\log c) ^2&
\log c  \log d\\ 
\log a  \log d &
\log b  \log d &
\log c  \log d &
(\log d) ^2
\end{pmatrix},
\end{align*}
whose eigenvalues and eigenvectors are
\begin{align*}
&0, 0, 0, 4 a^{2 x_0} b^{2 x_1} c^{2 x_2} d^{ 2 x_3} ((\log a)^2 + (\log b)^2 + (\log c)^2 + (\log d)^2), \\
&(-\log b, \log a, 0, 0), (-\log c, 0, \log a, 0),(-\log d, 0, 0, \log a), (\log a, \log b, \log c, \log d).
\end{align*}
Therefore $a^{x_0}b^{x_1}c^{x_2}d^{x_3}$ is convex since the eigenvalues are non-negative.

By the convexity of $T$ (in particular, the convexity with respect to $x_3$) we can write
\begin{align}
T_3(x_0,x_1,x_2)&:=\min_{\alpha_3 \in \mathbb R} \max_{\alpha_3 \leq x_3 \leq \alpha_3+1 } T(x_0,x_1,x_2,x_3) \notag \\
&=\min_{\alpha_3 \in \mathbb R}\max \{T(x_0,x_1,x_2,\alpha_3),T(x_0,x_1,x_2,\alpha_3+1)\} \label{eq1} \\
&=T(x_0,x_1,x_2,\alpha) \notag
\end{align}
for a unique $\alpha$ satisfying
\begin{align*}
T(x_0,x_1,x_2,\alpha) =T(x_0,x_1,x_2,\alpha+1).
\end{align*}
Now we clam that $T_3(x_0,x_1,x_2)$ is again convex: namely we have for $t \in [0,1]$
\begin{align*}
T_3(ta_0+(1-t)b_0,ta_1+(1-t)b_1,ta_2+(1-t)b_2) \leq tT_3(a_0,a_1,a_2)+(1-t)T_3(b_0,b_1,b_2).
\end{align*}
Say
\begin{align*}
&T_3(a_0,a_1,a_2)=T(a_0,a_1,a_2,\alpha)=T(a_0,a_1,a_2,\alpha+1), \\
&T_3(b_0,b_1,b_2)=T(b_0,b_1,b_2,\beta)=T(b_0,b_1,b_2,\beta+1).
\end{align*}
Moreover we put 
\begin{align*}
c_i:=ta_i+(1-t)b_i \ (i=0,1,2), \quad  c_3:=t\alpha+(1-t)\beta.
\end{align*}
Since 
\begin{align*}
&(a_0,a_1,a_2,\alpha,T_3(a_0,a_1,a_2)), &&(b_0,b_1,b_2,\beta,T_3(b_0,b_1,b_2)), \\
&(a_0,a_1,a_2,\alpha+1,T_3(a_0,a_1,a_2)), &&(b_0,b_1,b_2,\beta+1,T_3(b_0,b_1,b_2))
\end{align*}
are elements of the convex set $S$, so are
\begin{align*}
&(c_0,c_1,c_2,c_3,tT_3(a_0,a_1,a_2)+(1-t)T_3(b_0,b_1,b_2)), \\
&(c_0,c_1,c_2,c_3+1,tT_3(a_0,a_1,a_2)+(1-t)T_3(b_0,b_1,b_2)).
\end{align*}
Namely we have
\begin{align*}
tT_3(a_0,a_1,a_2)+(1-t)T_3(b_0,b_1,b_2)\geq T(c_0,c_1,c_2,c_3), T(c_0,c_1,c_2,c_3+1).
\end{align*}
Hence, by (\ref{eq1}), we have
\begin{align*}
tT_3(a_0,a_1,a_2)+(1-t)T_3(b_0,b_1,b_2)&\geq \max\{T(c_0,c_1,c_2,c_3), T(c_0,c_1,c_2,c_3+1)\} \\
&\geq T_3(c_0,c_1,c_2)
\end{align*}
as desired.
By repeating the same argument, we can write
\begin{align*}
&T_2(x_0,x_1):=\min_{\alpha_2 \in \mathbb R} \max_{\alpha_2 \leq x_2 \leq \alpha_2+1 } T_3(x_0,x_1,x_2)=T_3(x_0,x_1,\alpha) \\
&\text{for} \ \alpha \ \text{with} \ T_3(x_0,x_1,\alpha) =T_3(x_0,x_1,\alpha+1), \\[5pt]
&T_1(x_0):=\min_{\alpha_1 \in \mathbb R} \max_{\alpha_1 \leq x_1 \leq \alpha_1+1 } T_2(x_0,x_1)=T_2(x_0,\alpha') \\
&\text{for} \ \alpha' \ \text{with} \ T_2(x_0,\alpha') =T_2(x_0,\alpha'+1).
\end{align*}

We easily obtain an upper bound of such minimal values as follows.
Consider a closed-interval $I=[a,b]$ and divide it into $N+1$ pieces:
\begin{align*}
A:=\{a,a+(b-a)/N,a+2(b-a)/N,\dots,b-(b-a)/N,b\}.
\end{align*}
Then we see that 
\begin{align*}
T_3(x_0,x_1,x_2)&=\min_{\alpha_3 \in \mathbb R}\max \{T(x_0,x_1,x_2,\alpha_3),T(x_0,x_1,x_2,\alpha_3+1)\} \\
&\leq \min_{\alpha_3 \in A}\max \{T(x_0,x_1,x_2,\alpha_3),T(x_0,x_1,x_2,\alpha_3+1)\}.
\end{align*}
By repeating similar arguments, we obtain an upper bound of $L$ as
\begin{align*}
L\leq 
\min_{\alpha_0 \in A}\max_{x_0=\alpha_0,\alpha_0+1}
\min_{\alpha_1 \in A}\max_{x_1=\alpha_1,\alpha_1+1}
\min_{\alpha_2 \in A}\max_{x_2=\alpha_2,\alpha_2+1}
\min_{\alpha_3 \in A}\max_{x_3=\alpha_3,\alpha_3+1} 
T(x_0,x_1,x_2,x_3).
\end{align*}
Now, we put $[a,b]:=[\frac{-101}{100},\frac{99}{100}]$, $N:=32$.
Then numerically we have
\begin{align*}
&t(\tfrac{-404}{400})=887.4\dots, \ \cdots ,\ t(\tfrac{-229}{400})=312.9\dots,  t(\tfrac{-204}{400})=260.8\dots, t(\tfrac{-179}{400})=241.1\dots, \\
&\cdots ,\ t(\tfrac{171}{400})=239.1\dots,t(\tfrac{196}{400})=259.0\dots, t(\tfrac{221}{400})=308.8\dots,\ \cdots, \ t(\tfrac{396}{400})=1094.5\dots.
\end{align*}
where we put 
\begin{align*}
t(\alpha_0):=\min_{\alpha_1 \in A}\max_{x_1=\alpha_1,\alpha_1+1}
\min_{\alpha_2 \in A}\max_{x_2=\alpha_2,\alpha_2+1}
\min_{\alpha_3 \in A}\max_{x_3=\alpha_3,\alpha_3+1} 
T(\alpha_0,x_1,x_2,x_3).
\end{align*}
Hence we obtain $L\leq \max \{t(\frac{-204}{400}),t(\frac{196}{400})\}=260.8\dots$ as desired.
\end{proof}

\begin{rmk}
Summarizing the proof of {\rm\cref{n=3}}, we showed that there exists a fundamental domain $D$ of $\mathbb R^{2^{n-1}}$ modulo $\mathbb Z^{2^{n-1}}$ satisfying 
\begin{align*}
\max \left\{\sum_{i=0}^{2^n-1}\prod_{j=0}^{2^{n-1}-1} |\sigma^{i+j}(\e)|^{2x_j} \mid (x_i)_i \in D\right\} < 2^n(1+8c_n)
\end{align*}
for $n=3$, by considering the $\mathbb Z$-module structure of $RE_+^n$. 
When $n\geq 4$, it seems to have to consider its Galois module structure, not only the $\mathbb Z$-module structure, in order to studying the relation between {\rm\cref{cnj}} 
and the class number.
We provide some partial (and numerical) results in the proceeding sections.
\end{rmk}

\section{$l$-Indivisibility of $h_n$ by numerical calculations} \label{lind}

We give a demonstration of numerical checks of the $l$-indivisibility of $k_n$ for several $(l,n)$, by using \cref{cnj}.
More powerful results can be seen in \cite{H1,H2,FK1,FK2,FK3,MO1,MO2}, including (\ref{fkf}).
Let $l$ be an odd prime.
We put
\begin{align*}
A_n^\frac{1}{l}&:=\{x \in \mathbb R \mid x^l \in A_n\}.
\end{align*}
Since $A_n/\{\pm 1\}$ is a free abelian group generated by $\{\sigma^i(\e) \mid i=0,\dots,2^{n-1}-1\}$, 
we may identify the following three $\mathbb F_l[G_n]$-modules 
\begin{align*}
\begin{array}{ccccc}
\mathbb F_l[x]/(x^{2^{n-1}}+1) &\cong& A_n^\frac{1}{l}/A_n &\cong& \mathbb F_l^{2^n-1}, \\
\rotatebox[origin=c]{90}{$\in$} &&\rotatebox[origin=c]{90}{$\in$} &&\rotatebox[origin=c]{90}{$\in$}\\
\displaystyle \overline{\sum_{i=0}^{2^{n-1}-1}a_i x^i} & \leftrightarrow& \displaystyle \overline{\prod_{i=0}^{2^{n-1}-1}\sigma^i(\e)^{\frac{a_i}{l}}}& \leftrightarrow& (a_i)_{0\leq i \leq 2^{n-1}-1}.
\end{array}
\end{align*}
Here $\sigma$ acts on $\mathbb F_l^{2^n-1}$ by
\begin{align*}
\sigma (a_0,a_1,a_2,\dots,a_{2^n-2},a_{2^n-1})=(-a_{2^n-1},a_0,a_1,\dots,a_{2^n-3},a_{2^n-2}).
\end{align*}
$\mathbb F_l[G_n]$ acts on $\mathbb F_l[x]/(x^{2^{n-1}}+1)$ via
\begin{align*}
\mathbb F_l[G_n] \st{\sigma \mt x}\cong \mathbb F_l[x]/(x^{2^{n}}-1) \twoheadrightarrow \mathbb F_l[x]/(x^{2^{n-1}}+1),
\end{align*}
and hence we may also consider $\mathbb F_l[x]/(x^{2^{n-1}}+1),A_n^\frac{1}{l}/A_n,\mathbb F_l^{2^n-1}$ are $\mathbb F_l[x]$-modules where $x$ acts as $\sigma$.

By the Chinese remainder theorem, the irreducible decomposition of $\mathbb F_l[x]/(x^{2^{n-1}}+1)$ as a $\mathbb F_l[x]$-module is given as 
\begin{align*}
&\mathbb F_l[x]/(x^{2^{n-1}}+1)=\bigoplus_{f_i} M_f, \\
&M_{f_i}:=\tfrac{x^{2^{n-1}}+1}{f_i} \cdot \mathbb F_l[x]/(x^{2^{n-1}}+1) \ (\cong \mathbb F_l^{\deg f_i}),
\end{align*}
where $f_i$ runs over all irreducible polynomial $f_i \in \mathbb F_l[x]$ dividing $x^{2^{n-1}}+1$.
(Note that $x^{2^{n-1}}+1 \bmod l$ has no multiple roots.)
Taking a polynomial $g_{f_i} \in \mathbb F_l[x]$ satisfying $\tfrac{x^{2^{n-1}}+1}{f_i}\cdot g_{f_i} \equiv 1 \mod f_i$, the idempotent map is given explicitly as
\begin{align*}
e_{f_i} \colon \mathbb F_l[x]/(x^{2^{n-1}}+1) \twoheadrightarrow M_{f_i}, \quad  h \mt \tfrac{x^{2^{n-1}}+1}{f_i} g_{f_i} h.
\end{align*}

Now, we assume that $l \mid k_{n}$. It follows that there exists $\epsilon \in RE_n^+$ satisfying 
\begin{align*}
\epsilon \notin A_n,\quad \epsilon^l \in A_n.
\end{align*}
This element corresponds to a non-trivial element $\overline {g_\epsilon} \in  \mathbb F_l[x]/(x^{2^{n-1}}+1) \cong A_n^\frac{1}{l}/A_n$.
Then we can take $f_i$ so that $e_{f_i}(\overline {g_\epsilon} )\neq \overline 0$ since $\sum_i e_{f_i}(\overline {g_\epsilon} )=\overline {g_\epsilon}  \neq \overline 0$. 
For such $f_i$, the whole of $M_{f_i}$ is contained in $RE_n^+/A_n$,
since we have  
\begin{align} \label{cond}
\begin{array}{cccc}
\mathbb F_l[x]/(x^{2^{n-1}}+1) &\cong &A_n^\frac{1}{l}/A_n \\
\cup & & \cup \\
M_{f_i} & \cong & \left\{g \cdot \overline{\e} \in A_n^\frac{1}{l}/A_n \mid g \in M_{f_i}\right\} \\
\rotatebox[origin=c]{90}{$=$} && \rotatebox[origin=c]{90}{$=$} \\
\mathbb F_l[x] \cdot e_{f_i}(\overline {g_\epsilon} ) & \cong & \left(\tfrac{x^{2^{n-1}}+1}{f_i} g_{f_i}\mathbb F_l[x]\right) \cdot \overline{\epsilon} & \subset RE_n^+/A_n.
\end{array}
\end{align}
For $g=\sum_{i=0}^m a_i x^i \in \mathbb Z[x]$ (not only for elements $\in \mathbb F_l[x]$), we put
\begin{align*}
g \cdot \e:=\prod_{i=0}^m \sigma^i(\e)^{\frac{a_i}{l}} \in A_n^{\frac{1}{l}}.
\end{align*}
Then the following proposition follows form (\ref{cond}).

\begin{prp}
Assume that an odd prime $l$ divides $k_{n}$. Then there exists an irreducible polynomial $f \in \mathbb F_l[x]$ dividing $x^{2^{n-1}}+1$ satisfying
\begin{align*}
\left\{\overline{g \cdot \e} \in A_n^\frac{1}{l}/A_n  \;\middle|\; g \in \mathbb Z[x] \text{ with } g \bmod (x^{2^{n-1}}+1)\in M_f \right\} \subset RE_n^+.
\end{align*}
\end{prp}

We extend the trace map to 
\begin{align*}
\widetilde \Tr\big(\epsilon^{\frac{2}{l}}\big):=\sum_{i=0}^{2^n-1}\left(\sigma^i(\epsilon)^2\right)^{\frac{1}{l}} \quad (\epsilon \in A_n).
\end{align*}
By the above proposition, \cref{cnj} can be used for a numerical check of the indivisibility of the class numbers as follows.

\begin{thm} \label{k=1}
Assume that {\rm\cref{cnj}} holds true for $n$. 
If for each irreducible polynomial $f \in \mathbb F_l[x]$ dividing $x^{2^{n-1}}+1$ there exists $g \in \mathbb Z[x]$ satisfying  
\begin{align*}
&g \bmod (l,x^{2^{n-1}}+1) \in M_f-\{0\}, \\
&\widetilde \Tr \left(\left(g \cdot \e\right)^2\right) < 2^n(1+8c_n),
\end{align*}
then we have  $l \nmid k_{n}$.
\end{thm}

\subsection{The case $n= 4,5$} \label{ind}

\begin{exm} \label{n=4}
Let $n=4$, $l<10^6$. For each irreducible polynomial $f \in \mathbb F_l[x]$ dividing $x^{2^3}+1$, 
we took the center lift $g$ of a suitable element in $M_f$ and confirmed that
\begin{align} \label{cnd4}
\Tr \left(\left(g \cdot \varepsilon_4 \right)^2\right) < 2^4(1+8c_4)=784.
\end{align}
Namely, by {\rm\cref{k=1}}, we checked that {\rm\cref{cnj}} implies $l \nmid k_4$ for $l<10^6$.

For example, let $l=3$.
Then the irreducible decomposition of $x^8+1 \bmod 3$ is given by
\begin{align*}
x^8+1 \equiv f_1f_2 \mod 3, \quad  f_1=x^4 + x^2 - 1,\quad f_2=x^4 - x^2 - 1.
\end{align*}
We choose elements $\frac{x^8+1}{f_i} \in M_{f_i}$ $(i=1,2)$ and take their center lifts $g_1:=x^4 - x^2 - 1$, $g_2:=x^4 + x^2 - 1$.
Then, by numerical computation,  we obtain 
\begin{align*}
&\widetilde \Tr \left(\left(g_1 \cdot \varepsilon_4\right)^2\right)=\widetilde \Tr \left(\left(\varepsilon_4^{-1}\sigma^2(\varepsilon_4)^{-1}\sigma^4(\varepsilon_4)\right)^\frac{2}{3}\right)
=95.6\dots, \\
&\widetilde \Tr \left(\left(g_2 \cdot \varepsilon_4\right)^2\right)=\widetilde \Tr \left(\left(\varepsilon_4^{-1}\sigma^2(\varepsilon_4)\sigma^4(\varepsilon_4)\right)^\frac{2}{3}\right)
=100.1\dots.
\end{align*}
These values satisfy the condition {\rm(\ref{cnd4})} for $3 \nmid k_4$.

Next, let $l=7$. Then we have
\begin{align*}
&x^8+1 \equiv f_1f_2f_3f_4 \mod 7, \\
&f_1=x^2 + x - 1, \quad f_2=x^2 + 3x - 1, \quad f_3=x^2 - 3x - 1, \quad f_4=x^2 - x - 1.
\end{align*}
First we take center lifts $g_i$ of $\frac{x^8+1}{f_i} \bmod 7$:
\begin{align*}
&g_1=x^6 - x^5 + 2x^4 - 3x^3 - 2x^2 - x - 1, \quad g_2=x^6 - 3x^5 + 3x^4 + 2x^3 - 3x^2 - 3x - 1, \\
&g_3=x^6 + 3x^5 + 3x^4 - 2x^3 - 3x^2 + 3x - 1,  \quad g_4=x^6 + x^5 + 2x^4 + 3x^3 - 2x^2 + x - 1.
\end{align*}
Then we have
\begin{align*}
&\widetilde \Tr \left(\left(g_1 \cdot \varepsilon_4\right)^2\right)=106.5\dots, \quad 
\widetilde \Tr \left(\left(g_2 \cdot \varepsilon_4\right)^2\right)=546.9\dots, \\
&\widetilde \Tr \left(\left(g_3 \cdot \varepsilon_4\right)^2\right)=840.6\dots, \quad 
\widetilde \Tr \left(\left(g_4 \cdot \varepsilon_4\right)^2\right)=160.2\dots.
\end{align*}
Note that the case $i=3$ does not satisfy the condition {\rm(\ref{cnd4})}. Replacing $g_3$ with the center lift $g_3'=2x^6 -x^5 -x^4 +3x^3 +x^2 -x - 2$ of $2\cdot \frac{x^8+1}{f_i} \bmod 7$, 
we have
\begin{align*}
\widetilde \Tr \left(\left(g_3' \cdot \varepsilon_4\right)^2\right)=200.7\dots, 
\end{align*}
which implies $7 \nmid k_4$.
\end{exm}

\begin{exm}
Let $n=5$, $l<10^5$, $\neq 97,193,257$.
Then, similarly as in {\rm\cref{n=4}}, the center lift $g$ of a suitable element in $M_f$ for each $f$ satisfies the condition of {\rm\cref{k=1}}:
\begin{align*}
\widetilde\Tr \left(\left(g \cdot \varepsilon_5 \right)^2\right) < 2^5(1+8c_5)=3104.
\end{align*}
We also check the exceptions $97,193,257$ by taking certain non-center lifts.  

Let $l=97$. Then we have
\begin{align*}
&x^{2^4}+1 \equiv \prod_{i=1}^{16} f_i \mod 97, \hspace{-100pt} \\
&f_1=x + 19, &&f_2=x + 20, &&f_3=x + 28, &&f_4=x + 30, &&f_5=x + 34, \\
&f_6=x + 42, &&f_7=x + 45, &&f_8=x + 46, &&f_9=x - 46, &&f_{10}=x - 45, \\
&f_{11}=x - 42, &&f_{12}=x - 34, &&f_{13}=x - 30, &&f_{14}=x - 28, &&f_{15}=x - 20, \\
&f_{16}=x - 19.
\end{align*}
For $i=1,4,5,6,8,9,11,13,14$, we put $g_i$ to be the center lift of $4\cdot \frac{x^{2^4}+1}{f_i}$. 
Then we have
\begin{align*}
&\widetilde\Tr \left( \left(g_{1}\cdot\varepsilon_5\right)^2\right)=1123.9\dots, 
&&\widetilde\Tr \left( \left(g_{4}\cdot\varepsilon_5\right)^2\right)=1429.9\dots, 
&&\widetilde\Tr \left( \left(g_{5}\cdot\varepsilon_5\right)^2\right)=2421.7\dots, \\
&\widetilde\Tr \left( \left(g_{6}\cdot\varepsilon_5\right)^2\right)=1632.8\dots, 
&&\widetilde\Tr \left( \left(g_{8}\cdot\varepsilon_5\right)^2\right)=2332.6\dots, 
&&\widetilde\Tr \left( \left(g_{9}\cdot\varepsilon_5\right)^2\right)=1291.7\dots, \\
&\widetilde\Tr \left( \left(g_{11}\cdot\varepsilon_5\right)^2\right)=1537.1\dots, 
&&\widetilde\Tr \left( \left(g_{13}\cdot\varepsilon_5\right)^2\right)=1492.2\dots, 
&&\widetilde\Tr \left( \left(g_{14}\cdot\varepsilon_5\right)^2\right)=1444.4\dots.
\end{align*}
For $i=2,3,7,10,12,15,16$, we have to take non-center lifts.
Hereinafter, we write a polynomial $\sum_{i=0}^k a_ix^i$ as a vector $[a_0,\dots,a_k]$ for saving pages.
We put 
\begin{align*}
g_2&:=[34, 8, 19, 33, 42, 27, -45, -22, -28, -18, 30, \underline{-50}, -46, 12, -20, 1], \\
g_3&:=[7, 24, 13, 3, 38, \underline{61}, 29, 44, 40, -43, 5, 31, 37, 16, 41, 2], \\
g_7&:=[41, -16, 37, -31, 5, 43, 40, \underline{53}, 29, 36, 38, -3, 13, -24, 7, 2], \\
g_{10}&:=[-41, -16, -37, -31, -5, \underline{-54}, -40, -44, -29, 36, -38, -3, -13, -24, -7, 2], \\
g_{12}&:=[17, \underline{49}, 10, 6, 23, -25, -15, -9, 14, -11, -26, -35, -21, -32, 39, 4], \\
g_{15}&:=[29, 16, -38, -31, \underline{-84}, -43, -7, -44, -41, -36, 37, -3, -5, 24, 40, 2], \\
g_{16}&:=[5, 36, 7, 31, -29, 24, 37, \underline{53}, 13, 16, -40, 3, 41, 43, 38, 2],
\end{align*}
which are lifts of 
$\frac{x^{2^4}+1}{f_2},2\cdot\frac{x^{2^4}+1}{f_3},2\cdot\frac{x^{2^4}+1}{f_7},2\cdot\frac{x^{2^4}+1}{f_{10}},4\cdot\frac{x^{2^4}+1}{f_{12}},2\cdot\frac{x^{2^4}+1}{f_{15}},2\cdot\frac{x^{2^4}+1}{f_{16}}$ 
respectively.
Here components with underlining are not contained in $[-\frac{l-1}{2},\frac{l-1}{2}]$.
Then we have
\begin{align*}
&\widetilde\Tr \left( \left(g_{2}\cdot\varepsilon_5\right)^2\right)=1492.1\dots, 
&&\widetilde\Tr \left( \left(g_{3}\cdot\varepsilon_5\right)^2\right)=1963.0\dots, 
&&\widetilde\Tr \left( \left(g_{7}\cdot\varepsilon_5\right)^2\right)=1548.9\dots, \\
&\widetilde\Tr \left( \left(g_{10}\cdot\varepsilon_5\right)^2\right)=920.6\dots, 
&&\widetilde\Tr \left( \left(g_{12}\cdot\varepsilon_5\right)^2\right)=1831.2\dots, 
&&\widetilde\Tr \left( \left(g_{15}\cdot\varepsilon_5\right)^2\right)=2985.0\dots, \\
&\widetilde\Tr \left( \left(g_{16}\cdot\varepsilon_5\right)^2\right)=2386.1\dots.
\end{align*}
The other cases $l=193,257$ can be done similarly.
\end{exm}

\begin{rmk}
Let $n=6$. Then, for many $l$ (e.g., $l=31,97,127,193,223,257,449,\dots$), the center lift $g$ of any element in $M_f$ does not satisfy the condition
\begin{align*}
\widetilde\Tr \left(\left(g \cdot \varepsilon_6 \right)^2\right) < 2^6(1+8c_6)=13376.
\end{align*} 
Moreover searching all non-center lifts is difficult due to the high dimension.
We confirmed that {\rm\cref{cnj}} implies that $l \nmid k_6$ only for $l=31$.
\end{rmk}

\subsection{The case $n=7$, $l>10^9$, $l\equiv 65 \bmod 128$} \label{new}

If $n,l$ are large, it is difficult to check the condition in \cref{k=1}.
However that becomes relatively easy in some special cases.
Let $n=7$, $l\equiv 65 \bmod 128$. 
We note that such $l$ are out of the range of (\ref{previous}).
Then the irreducible decomposition of $x^{2^{6}} +1 \bmod l$ is in the form
\begin{align} \label{dcm}
x^{2^{6}} +1 \bmod l=\prod_{i=1}^{32} (x^2+a_i).
\end{align}
In fact, that $l\equiv 65 \bmod 128$ is equivalent to that $l$ splits completely in $\mathbb Q(\zeta_6)$ and does not in $\mathbb Q(\zeta_7)$ where $\zeta_{n}:=e^{\frac{2\pi i}{2^n}}$.
Then $y^{2^5} +1$, which is a minimal polynomial of $\mathbb Q(\zeta_6)$, decomposes a product of polynomials of degree $1$ modulo $l$, 
and $x^{2^6} +1$ of $\mathbb Q(\zeta_7)$ does not. Considering $y=x^2$, we obtain the expression (\ref{dcm}).
Since half of the coefficients of $\frac{x^{2^{6}} +1}{x^2+a}$ are equal to $0$ (that is, $\frac{x^{2^{6}} +1}{x^2+a}$ is in the form $\sum_{i=0}^{31}c_{2i} x^{2i}$), 
the value of $\widetilde\Tr ( (g\cdot\varepsilon_7)^2)$ tends to ``small'' if we take a center lift $g$ of $\frac{x^{2^{6}} +1}{x^2+a}$ multiplied by a constant.
For example, let $l=1000000321$, which is the least prime satisfying $l>10^9$, $l\equiv 65 \bmod 128$. Note that this case is not contained in (\ref{previous}).
Then we have 
\begin{align*}
&x^{2^{6}} +1 \bmod l = \prod_{i=1}^{32} (x^2+a_i), \\
&a_1=30063488, \ a_2=30912022, \ a_3=42483948, \ a_4=59955883, \ a_5=78186285, \\
&a_6=160612070, \ a_7=191346380, \ a_8=246360387, \ a_9=268629094, \ a_{10}=269645956, \\
&a_{11}=280492327, \ a_{12}=303644312, \ a_{13}=311722386, \ a_{14}=424439170, \\
&a_{15}=441230693, \ a_{16}=447503416, \ a_{16+i}=-a_{17-i} \ (1\leq i \leq 16).
\end{align*}
We put $g_i$ to be the center lift of $b_i\cdot \frac{x^{2^{6}} +1}{x^2+a_i}$ with 
\begin{align*}
&b_{1}=231,\ b_{2}=231,\ b_{3}=867,\ b_{4}=125,\ b_{5}=386,\ b_{6}=231,\ b_{7}=100,\ b_{8}=100,\\
&b_{9}=64, \ \ b_{10}=36,\ b_{11}=702,\ b_{12}=771,\ b_{13}=231,\ b_{14}=2069,\ b_{15}=349,\ b_{16}=64,\\
&b_{17}=64,\ b_{18}=64,\ b_{19}=4,\ b_{20}=64,\ b_{21}=686,\ b_{22}=105,\ b_{23}=167,\ b_{24}=64,\\
&b_{25}=100,\ b_{26}=89,\ b_{27}=100,\ b_{28}=100,\ b_{29}=100,\ b_{30}=100,\ b_{31}=100,\ b_{32}=64.
\end{align*}
Then $t_i:=\widetilde\Tr ((g_{i}\cdot\varepsilon_7)^2)$ are calculated numerically as follows.
\begin{align*}
&t_{1}=24947.7\dots,\ 
t_{2}=15616.7\dots,\ 
t_{3}=49165.2\dots,\ 
t_{4}=23454.0\dots,\ 
t_{5}=46028.1\dots,\\ 
&t_{6}=41400.4\dots,\ 
t_{7}=19344.5\dots,\ 
t_{8}=26943.5\dots,\ 
t_{9}=42868.4\dots,\ 
t_{10}=40913.4\dots,\\ 
&t_{11}=44067.7\dots,\ 
t_{12}=49457.9\dots,\ 
t_{13}=18759.3\dots,\ 
t_{14}=39188.3\dots,\\ 
&t_{15}=35939.1\dots,\ 
t_{16}=44713.3\dots,\ 
t_{17}=41782.1\dots,\ 
t_{18}=47974.8\dots,\\ 
&t_{19}=52445.8\dots,\ 
t_{20}=49841.0\dots,\ 
t_{21}=43256.3\dots,\ 
t_{22}=52244.6\dots,\\ 
&t_{23}=49338.6\dots,\ 
t_{24}=22229.3\dots,\ 
t_{25}=36290.0\dots,\ 
t_{26}=48593.0\dots,\\ 
&t_{27}=26438.3\dots,\ 
t_{28}=40208.3\dots,\ 
t_{29}=23006.2\dots,\ 
t_{30}=19831.0\dots,\\ 
&t_{31}=16060.6\dots,\ 
t_{32}=42470.9\dots.
\end{align*}
These values satisfy the condition in \cref{k=1}:
\begin{align*}
\widetilde\Tr \left(\left(g \cdot \varepsilon_7 \right)^2\right) < 2^7(1+8c_7)=53376.
\end{align*} 
Namely, \cref{cnj} for $n=7$ implies $l=1000000321 \nmid k_7$. 
Similarly, we checked that \cref{cnj} for $n=7$ implies the $l$-indivisibility of $k_7$ for first $1000$ primes satisfying 
\begin{align*}
10^9<l,\ l\equiv 65 \bmod 128,
\end{align*}
form $1000000321$ to $1001287361$.

\end{document}